\documentclass{amsart}
\usepackage{amsfonts}

\usepackage{graphicx}
\usepackage{amscd}
\usepackage{amsmath}

\newtheorem{theorem}{Theorem}
\theoremstyle{plain}

\newtheorem{corollary}{Corollary}

\newtheorem{lemma}{Lemma}

\newtheorem{proposition}{Proposition}

\numberwithin{equation}{section}

\begin{document}
\author{}
\thanks{Universit\`{a} degli Studi di Bergamo, Facolt\`{a} di Ingegneria, \noindent
Viale Marconi, 5, 24044 Dalmine \noindent(Bergamo) Italy }
\thanks{E-mail: Garattini@mi.infn.it }
\title{}
\maketitle

\begin{center}
\bigskip {\large {\textbf{REMO GARATTINI}}} \vspace{1cm}

\textrm{\textbf{{\large {On the existence of nontrivial solutions for a
nonlinear equation relative to a measure-valued Lagrangian on homogeneous
spaces}}}}\vspace{1cm}
\end{center}

\noindent \textbf{Abstract.} We prove the existence of a non-trivial
solution for a nonlinear equation related to a measure-valued Lagrangian.
The result is based on a compact embedding theorem of the Lagrangian domain
and on the application of the Mountain Pass Theorem joined to a Palais-Smale
condition.

\section{Introduction and Result}

We consider a locally compact separable Hausdorff topological space $X$
endowed with a measure $m$ and a quasidistance $d$. A quasidistance $d$ on $%
X $ is a function on $X\times X$ with the usual properties of a metric and a
weaker version of the triangle inequality
\begin{equation*}
d\left( x,y\right) \leq c_{T}\left( d\left( \left( x,z\right) +d\left(
z,y\right) \right) \right) ,\qquad c_{T}\geq 1.
\end{equation*}
The set
\begin{equation*}
B\left( x,R\right) =\left\{ y\in X:d\left( x,y\right) <R\right\}
\end{equation*}
will be called a quasi-ball. The triple $\left( X,d,m\right) $ is assumed to
satisfy the following property: for every $R_{0}>0$ there exists a constant $%
c_{0}>0,$ dependent on $R_{0}$, such that for $r\leq \frac{R}{2}\leq R\leq
R_{0}$%
\begin{equation}
0<c_{0}\left( \frac{r}{R}\right) ^{\nu }m\left( B\left( x,R\right) \right)
\leq m\left( B\left( x,r\right) \right)  \label{a1}
\end{equation}
for every $x\in X$, where $\nu $ is a positive real number independent of $%
r,R,R_{0}$. Such a triple $\left( X,d,m\right) $ will be called a
homogeneous space of dimension $\nu .$ We point out, however, that a given
exponent $\nu $ occurring in $\left( \ref{a1}\right) $ should be considered,
more precisely, as an upper bound of the ``homogeneous dimension'', hence we
should better call $\left( X,d,m\right) $ a homogeneous space of dimension
less or equal than $\nu .$ Our setting is given by a couple $\left( X,%
\mathcal{L}\right) $, ``a homogeneous space $X$ with a Lagrangian $\mathcal{L%
}$'', with the following properties

\begin{description}
\item[(L1)]  $\mathcal{L}:\mathcal{C}\longmapsto \mathcal{M}\left( X\right) $
is a map which associates with each function $u$ from a given subspace $%
\mathcal{C}$ of $C\left( X\right) $ a measure $\mathcal{L}\left[ u\right]
\in \mathcal{M}^{+}\left( X\right) $, where $C\left( X\right) $ denotes the
space of all continuous functions on $X$ and $\mathcal{M}^{+}\left( X\right)
$ the space of all nonnegative Radon measures on $X$.

\item[(L2)]  We assume that there exists $k\geq 1$ such that for a given $%
p\geq 1$, the following family of Poincar\'{e}-like inequalities holds on
the metric quasi-balls $B\left( x,r\right) \subset \subset X$ \cite
{Biroli-Mosco}\cite{Maly-Mosco}:
\begin{equation}
\int\limits_{B\left( x,r\right) }\left| u-u_{x,r}\right| ^{p}dm\leq
c_{P}r^{p}\int\limits_{B\left( x,kr\right) }d\mathcal{L}\left[ u\right] ,
\label{a2}
\end{equation}
where $u_{x,r}$ is the average of $u$ on $B\left( x,r\right) $, for every $%
u\in \mathcal{C}$ and $B\left( x,r\right) \subset \subset X$.

\item[(L3)]  If $u\in \mathcal{C}$ and $g\in C^{1}\left( \mathbf{R}\right) $
with $g^{\prime }$ bounded on $\mathbf{R}$, then $g\left( u\right)
:x\longmapsto g\left( u\left( x\right) \right) $ belong also to $\mathcal{C}$
and
\begin{equation}
\mathcal{L}\left[ g\left( u\right) \right] =\left| g^{\prime }\left(
u\right) \right| ^{p}\mathcal{L}\left[ u\right]
\end{equation}
\end{description}

We are interested in nontrivial solution of the following problem
\begin{equation}
\int\limits_{X}d\mathcal{L}\left[ u\right] v\left( x\right)
+\int\limits_{X}V\left( x\right) u^{p}\left( x\right) v\left( x\right)
m\left( dx\right) =\int\limits_{X}f\left( u\left( x\right) \right) v\left(
x\right) m\left( dx\right)   \label{a3}
\end{equation}
for every $v\in \mathcal{C}\cap L^{p}\left( X,Vm\right) $ where $u\in
\mathcal{C}\cap L^{p}\left( X,Vm\right) $ ( $Vm$ is the Radon measure with
density $V$ with respect to $m$). Eq. $\left( \ref{a3}\right) $ is a
generalization of the problem of searching for nontrivial solution for a
semilinear equation in the framework of Dirichlet forms as studied in Ref.
\cite{Biroli-Tersian} and in the framework of semilinear equations of the
form
\begin{equation}
\triangle u+u^{p}=0
\end{equation}
considered in Ref.\cite{Falconer}. Further developments on semilinear
equations for Dirichlet forms can be found in Ref.\cite{Matzeu} for problems
of the type
\begin{equation*}
\int\limits_{\Omega }\alpha \left( u,v\right) \left( dx\right) -\lambda
\int\limits_{\Omega }a\left( x\right) u\left( x\right) v\left( x\right)
m\left( dx\right) =\int\limits_{\Omega }f\left( u\left( x\right) \right)
v\left( x\right) m\left( dx\right) ,
\end{equation*}
where $\Omega $ is an open bounded subset of $X$, $\alpha \left( u,v\right) $
is a uniquely defined signed Radon measure on $X$, $\lambda $ is an
arbitrary nonvanishing number and $a\in Lip\left( \bar{\Omega}\right) $ with
$a\left( x\right) >0$. To analyze Eq. $\left( \ref{a3}\right) $, we assume
that
\begin{equation}
W=\left\{ u:\int\limits_{X}d\mathcal{L}\left[ u\right] +\int%
\limits_{X}Vu^{p}m\left( dx\right) <+\infty \right\}   \label{a3a}
\end{equation}
and that
\begin{equation}
\left\| u\right\| _{W}=\left[ \int\limits_{X}d\mathcal{L}\left[ u\right]
+\int\limits_{X}Vu^{p}m\left( dx\right) \right] ^{\frac{1}{p}}
\end{equation}
be a norm in $W$. Moreover let us assume that $V\in C\left( X,\mathbb{R}%
\right) $ and
\begin{equation}
V\left( x\right) >0,\qquad \forall x\in X  \label{a4}
\end{equation}
\begin{equation}
V\left( x\right) \rightarrow +\infty ,\qquad \text{as\qquad }d\left(
0,x\right) \rightarrow +\infty   \label{a5}
\end{equation}
where $0$ is an arbitrarily fixed point in $X$. We assume also that $f\left(
t\right) \in C\left( X,\mathbb{R}\right) $ satisfies the following conditions
\begin{equation}
f\left( 0\right) =0,\ f\left( t\right) =o\left( t\right) ,\qquad \text{as }%
t\rightarrow 0  \label{a6}
\end{equation}
\begin{equation}
f\left( t\right) =o\left( \left| t\right| ^{\frac{\nu +p}{\nu -p}}\right)
,\qquad \text{as }\left| t\right| \rightarrow +\infty   \label{a7}
\end{equation}
if $\nu >p$ or
\begin{equation}
f\left( t\right) =o\left( \left| t\right| ^{\sigma }\right) ,\qquad \text{as
}\left| t\right| \rightarrow +\infty   \label{a8}
\end{equation}
$\sigma >p+1$, if $\nu \leq p$. Finally we assume that
\begin{equation}
0<\mu F\left( t\right) =\mu \int_{0}^{t}f\left( s\right) ds\leq tf\left(
t\right)   \label{a9}
\end{equation}
where $p<\frac{p\nu }{\nu -p}$ if $\nu >p$ or $p<\mu $ if $\nu \leq p$. We
observe that from the assumption $\left( \ref{a9}\right) $ it follows that
there exists $m>0$ such that
\begin{equation}
F\left( t\right) \geq m\left| t\right| ^{\mu }
\end{equation}
for $\left| t\right| \geq 1$. The result we will prove in the next Section
is the following:

\begin{theorem}
\label{t1}Let the assumptions $\left( \ref{a4}\right) ,\left( \ref{a5}%
\right) ,\left( \ref{a6}\right) ,\left( \ref{a9}\right) $ hold together with
$\left( \ref{a7}\right) $ if $\nu >2$ or with $\left( \ref{a8}\right) $ if $%
\nu =2$. Then the problem $\left( \ref{a3}\right) $ has a nontrivial
solution.
\end{theorem}

\textbf{Acknowledgments:} The Author wishes to thank Marco Biroli for very
useful and stimulating discussions on this subject.

\section{Preliminary results}

We begin the section with a covering Lemma and its Corollary.

\begin{lemma}
\label{l0}A ball $B\left( x,R\right) $ can be covered by a finite number $%
n\left( r,R\right) $ of balls $B\left( x_{i},r\right) ,$ $r\leq R$, such
that $x_{i}\in B\left( x,R\right) $ and $B\left( x_{i},\frac{r}{2}\right)
\cap B\left( x_{j},\frac{r}{2}\right) =\emptyset $ for $i\neq j$. Moreover
every point of $B\left( x,R\right) $ is covered by at most $M$ balls $%
B\left( x_{i},R\right) $ where $M$ depends on $r$.
\end{lemma}

\begin{proof}
The first part of the result follows immediately from assumption $\left( \ref
{a1}\right) $. For the second part we observe that if a point $x$ in $%
B\left( x,R\right) $ is covered by the ball $B\left( x_{i},r\right) $, then $%
x_{i}\in B\left( x,r\right) $; so the number $M$ of the balls $B\left(
x_{i},r\right) $, that cover $x$, is estimated by the greatest number $Q$ of
points $y_{k}$ in $B\left( x,r\right) $ with $d\left(
y_{k_{1}},y_{k_{2}}\right) \geq \frac{r}{2}$ and we observe that, by $\left(
\ref{a1}\right) $, $Q$ is estimated by a number $M$ depending only on $r$.
\end{proof}

From Lemma \ref{l0}, we obtain the following

\begin{corollary}
The space $X$ can be covered by a countable union of balls $B\left(
x_{i},r\right) $, such that $B\left( x_{i},\frac{r}{2}\right) \cap B\left(
x_{j},\frac{r}{2}\right) =\emptyset $ for $i\neq j$. Moreover every point of
$X$ is covered by at most $M$ balls, where $M$ depends only on $r$.
\end{corollary}

We prove now a compact embedding result

\begin{lemma}
\label{l1}Let the assumption related to inequality $\left( \ref{a2}\right) $
holds. Then every sequence $\left\{ u_{n}\right\} $ in $\mathcal{C}$ $\left[
B\left( x,\left( k+1\right) R\right) \right] $ such that
\begin{equation}
\int\limits_{B\left( x,kr\right) }d\mathcal{L}\left[ u\right] \leq C
\end{equation}
is relatively compact in $L^{p}\left( B\left( x,R\right) ,m\right) .$
\end{lemma}

\begin{proof}
We have to prove that there is a subsequence of $\left\{ u_{n}\right\} $
convergent in $L^{p}\left( B\left( x,R\right) ,m\right) $. Taking into
account assumption $\left( \ref{a1}\right) $, the ball $B\left( x,R\right) $
can be covered by a finite number of balls $B\left( x_{i},r\right) $, $r\leq
\frac{R}{4}$, $j=1,\ldots ,Q$ where $Q$ depends on $r,R$, such that every
point of $B\left( x,R\right) $ belongs at most to $M$ balls, where $M$ does
not depend on $r$. Let $w_{n,m}=u_{n}-u_{m}$ and $\bar{w}_{n,m}=%
\int\limits_{B\left(x_{j},r\right) }\hspace{-0.31cm}%
%
w_{n,m}m\left( dx\right) .$ Then
\begin{equation*}
\int\limits_{B\left( x,R\right) }w_{n,m}^{p}m\left( dx\right) \leq
\sum_{j=1}^{Q}\int\limits_{B\left( x_{j},r\right) }w_{n,m}^{p}m\left(
dx\right) =\sum_{j=1}^{Q}\int\limits_{B\left( x_{j},r\right) }\left| w_{n,m}-%
\bar{w}_{n,m}+\bar{w}_{n,m}\right| ^{p}m\left( dx\right)
\end{equation*}
\begin{equation}
\leq 2^{p-1}\sum_{j=1}^{Q}\int\limits_{B\left( x_{j},r\right) }\left|
w_{n,m}-\bar{w}_{n,m}\right| ^{p}m\left( dx\right)
+2^{p-1}\sum_{j=1}^{Q}\int\limits_{B\left( x_{j},r\right) }\left( \bar{w}%
_{n,m}\right) ^{p}m\left( dx\right) .  \label{1}
\end{equation}
Since
\begin{equation*}
\int\limits_{B\left( x_{j},r\right) }\left( \bar{w}_{n,m}\right) ^{p}m\left(
dx\right) =\int\limits_{B\left( x_{j},r\right) }\frac{m\left( dx\right) }{%
m^{p}\left( B\left( x_{j},r\right) \right) }\left( \int\limits_{B\left(
x_{j},r\right) }\left( w_{n,m}\right) m\left( dx\right) \right) ^{p}
\end{equation*}
\begin{equation}
=\frac{1}{m^{p-1}\left( B\left( x_{j},r\right) \right) }\left(
\int\limits_{B\left( x_{j},r\right) }\left( w_{n,m}\right) m\left( dx\right)
\right) ^{p},
\end{equation}
then inequality $\left( \ref{1}\right) $ becomes
\begin{equation*}
2^{p-1}\sum_{j=1}^{Q}\int\limits_{B\left( x_{j},r\right) }\left| w_{n,m}-%
\bar{w}_{n,m}\right| ^{p}m\left( dx\right)
+2^{p-1}\sum_{j=1}^{Q}\int\limits_{B\left( x_{j},r\right) }\left( \bar{w}%
_{n,m}\right) ^{p}m\left( dx\right)
\end{equation*}
\begin{equation*}
\leq 2^{p-1}c_{p}r^{\alpha }\sum_{j=1}^{Q}\int\limits_{B\left(
x_{j},kr\right) }d\mathcal{L}\left[ u\right] +2^{p-1}\sum_{j=1}^{Q}\frac{1}{%
m^{p-1}\left( B\left( x_{j},r\right) \right) }\left( \int\limits_{B\left(
x_{j},r\right) }\left( w_{n,m}\right) m\left( dx\right) \right) ^{p}
\end{equation*}
\begin{equation}
\leq 2^{p-1}c_{p}r^{\alpha }MCk^{\nu }+\left( \frac{R}{r}\right) ^{\nu
\left( p-1\right) }\frac{2^{p-1}}{m^{p-1}\left( B\left( x,R\right) \right)
c_{0}}\sum_{j=1}^{Q}\left( \int\limits_{B\left( x_{j},r\right) }\left(
w_{n,m}\right) m\left( dx\right) \right) ^{p}.
\end{equation}
Choose $r=r_{\varepsilon }$ and $\varepsilon >0$ such that $%
2^{p-1}c_{p}r_{\varepsilon }^{\alpha }MCk^{\nu }\leq \frac{\varepsilon }{2}$%
. Suppose $\left\{ u_{n}\right\} $ is weakly convergent in $L^{p}\left(
B\left( x,\left( k+1\right) R\right) ,m\right) $ then
\begin{equation}
\left( \frac{R}{r_{\varepsilon }}\right) ^{\nu \left( p-1\right) }\frac{%
2^{p-1}}{m^{p-1}\left( B\left( x,R\right) \right) c_{0}}\sum_{j=1}^{Q}\left(
\int\limits_{B\left( x_{j},r\right) }\left( w_{n,m}\right) m\left( dx\right)
\right) ^{p}\leq \frac{\varepsilon }{2}
\end{equation}
for $n,m\geq n_{\varepsilon }$. This implies
\begin{equation}
\int\limits_{B\left( x,R\right) }w_{n,m}^{p}m\left( dx\right) \leq
\varepsilon
\end{equation}
and $\left\{ u_{n}\right\} $is a Cauchy sequence in the space $L^{p}\left(
B\left( x,R\right) ,m\right) $ then $\left\{ u_{n}\right\} $is convergent in
$L^{p}\left( B\left( x,R\right) ,m\right) $.
\end{proof}

\begin{lemma}
\label{l2}Let $W$ $\subset \mathcal{C}$ be the space defined in Eq. $\left(
\ref{a3a}\right) $ and let us assume that $W$ be a Banach space w.r.t. $%
\left\| .\right\| _{W}$, then the embedding of $W$ in $L^{p}\left(
X,m\right) $ is compact.
\end{lemma}

\begin{proof}
Let $\left\| u_{k}\right\| _{W}\leq C$. After extraction of a subsequence,
we have that $\left\{ u_{k}\right\} $ is weakly convergent in $W$ to $u$. We
suppose, without loss of generality that $u=0$ and prove
\begin{equation}
\int\limits_{X}u_{k}^{p}m\left( dx\right) \rightarrow 0
\end{equation}
when $k\rightarrow +\infty $. Let $\varepsilon >0$, $\exists $ $R>0$ such
that $V\left( x\right) \geq \frac{1+C^{p}}{\varepsilon }$ when $d\left(
x,0\right) \geq R$. Since $\int\limits_{B\left( 0,R\right) }u_{k}^{p}m\left(
dx\right) \rightarrow 0$ when $k\rightarrow +\infty $, then $\exists $ $k$
such that for $k\geq k_{\varepsilon }$%
\begin{equation}
\int\limits_{B\left( 0,R\right) }u_{k}^{p}m\left( dx\right) \leq \frac{%
\varepsilon }{1+C^{p}}.
\end{equation}
Then for $k\geq k_{\varepsilon }$%
\begin{equation*}
\int\limits_{X}u_{k}^{p}m\left( dx\right) =\int\limits_{B\left( 0,R\right)
}u_{k}^{p}m\left( dx\right) +\int\limits_{X\backslash B\left( 0,R\right)
}u_{k}^{p}m\left( dx\right)
\end{equation*}
\begin{equation*}
\leq \frac{\varepsilon }{1+C^{p}}+\int\limits_{X\backslash B\left(
0,R\right) }u_{k}^{p}m\left( dx\right) \leq \frac{\varepsilon }{1+C^{p}}%
\left[ 1+\int\limits_{X\backslash B\left( 0,R\right) }Vu_{k}^{p}m\left(
dx\right) \right]
\end{equation*}
\begin{equation}
\leq \frac{\varepsilon }{1+C^{p}}\left[ 1+\left\| u_{k}\right\| _{W}^{p}%
\right] \leq \varepsilon .
\end{equation}
\end{proof}

\section{Proof of Theorem\ref{t1}}

The function on $W$ associated to our problem can be written as
\begin{equation}
\varphi\left( u\right) =\frac{1}{2}\left\| u\right\|
_{W}^{p}-\int\limits_{X}F\left( u\left( x\right) \right) m\left( dx\right) .
\end{equation}
It can be proved that $\varphi\in C^{1}\left( W,\mathbb{R}\right) $ and
\begin{equation}
\left\langle \varphi^{\prime}\left( u\right) ,v\right\rangle =\left(
u,v\right) _{W}-\int\limits_{X}f\left( x,u\left( x\right) \right) v\left(
x\right) m\left( dx\right) .
\end{equation}
The critical points of $\varphi$ are weak solution of our problem, then to
prove Theorem\ref{t1} it is enough to prove the existence of nontrivial
points for $\varphi$.

\begin{proposition}
The functional $\varphi $ satisfies the Palais-Smale condition under
assumption of Theorem \ref{t1}
\end{proposition}

\begin{proof}
Let $\left\{ u_{k}\right\} $ be a sequence in $W$ such that
\begin{equation}
\left| \varphi \left( u_{k}\right) \right| \leq C\qquad \varphi ^{\prime
}\left( u_{k}\right) \rightarrow 0,  \label{p2}
\end{equation}
in $W^{\ast }$ as $k\rightarrow +\infty $, where $W^{\ast }$ denotes the
dual space of $W$. From $\left( \ref{p2}\right) $ we obtain that there
exists $k_{0}$ such that for $k\geq k_{0}$%
\begin{equation}
\left| \left\langle \varphi ^{\prime }\left( u_{k}\right)
,u_{k}\right\rangle \right| \leq \mu \left\| u_{k}\right\| _{W}.
\end{equation}
Then
\begin{equation*}
C+\left\| u\right\| _{W}^{p}\geq \varphi \left( u_{k}\right) -\frac{1}{\mu }%
\left\langle \varphi ^{\prime }\left( u_{k}\right) ,u_{k}\right\rangle
\end{equation*}
\begin{equation*}
=\frac{1}{2}\left\| u_{k}\right\| _{W}^{p}-\int\limits_{X}F\left(
u_{k}\left( x\right) \right) m\left( dx\right) -\frac{1}{\mu }\left( \left\|
u_{k}\right\| _{W}^{p}-\int\limits_{X}f\left( u_{k}\left( x\right) \right)
u_{k}m\left( dx\right) \right)
\end{equation*}
\begin{equation}
=\left( \frac{1}{2}-\frac{1}{\mu }\right) \left\| u_{k}\right\|
_{W}^{p}-\int\limits_{X}F\left( u_{k}\left( x\right) \right) m\left(
dx\right) -\frac{1}{\mu }\int\limits_{X}f\left( u_{k}\left( x\right) \right)
u_{k}m\left( dx\right) \geq \left( \frac{1}{2}-\frac{1}{\mu }\right) \left\|
u_{k}\right\| _{W}^{p}.
\end{equation}
Since $\left\{ u_{k}\right\} $ is bounded in $W$ and from Lemma \ref{l2}, we
know that there exists a subsequence strongly convergent in $L^{p}\left(
X,m\right) $ and weakly to $u\in W.$ We apply now the Lemma 5 if $\nu \geq p$
or the Lemma 6 if $\nu <p$ of Ref.\cite{Biroli-Tersian} to the function $%
g\left( t\right) =f\left( t\right) $ and to the sequence $\left(
u_{k}\right) $ and we obtain
\begin{equation}
\lim_{k\rightarrow +\infty }\int\limits_{X}f\left( u_{k}\right) \left(
u_{k}-u\right) m\left( dx\right) =0.  \label{3}
\end{equation}
From the assumption we have that
\begin{equation}
\left| \left\langle \varphi ^{\prime }\left( u\right) ,v\right\rangle
\right| \leq \varepsilon _{k}\left\| v\right\| _{W}  \label{4}
\end{equation}
where $\varepsilon _{k}\rightarrow 0$ as $k\rightarrow +\infty $. Then from $%
\left( \ref{4}\right) $ we have
\begin{equation*}
\left\langle \varphi ^{\prime }\left( u_{k}\right) ,u_{k}-u\right\rangle
=\left( u_{k},u_{k}-u\right) _{W}-\int\limits_{X}f\left( x,u_{k}\left(
x\right) \right) \left( u_{k}-u\right) \left( x\right) m\left( dx\right)
\end{equation*}
\begin{equation}
=\left\| u_{k}\right\| _{W}^{p}-\left( u_{k},u\right)
_{W}-\int\limits_{X}f\left( x,u_{k}\left( x\right) \right) \left(
u_{k}-u\right) \left( x\right) m\left( dx\right) .
\end{equation}
From $\left( \ref{3}\right) $ and $\left( \ref{4}\right) $ we obtain
\begin{equation}
\left\langle \varphi ^{\prime }\left( u_{k}\right) ,u_{k}-u\right\rangle
\rightarrow \left\| u_{k}\right\| _{W}^{p}-\left( u_{k},u\right)
_{W}\rightarrow 0,
\end{equation}
when $k\rightarrow +\infty $. This implies that $\left\{ u_{k}\right\} $
converges to $u$ strongly in $W$.
\end{proof}

\begin{proof}[Proof of Theorem\ref{t1}]
First we prove that for $\rho \leq \min \left( \frac{a}{2}m\left( B\left(
0,1\right) \right) ,\frac{1}{2}\right) $ small enough $\varphi \left(
u\right) \geq \gamma >0$ for $\left\| u_{k}\right\| _{W}=\rho $. Consider
the case $\nu \geq p$. As in Lemma 5 of Ref.\cite{Biroli-Tersian} we obtain
that for every $\varepsilon >0$ there exists a constant $C_{\varepsilon }$
such that
\begin{equation}
0\leq F\left( t\right) \leq \varepsilon \left( \left| t\right| ^{p}+\left|
t\right| ^{\beta }\right) +C_{\varepsilon }\left| t\right| ^{\beta }
\end{equation}
where $\beta =\frac{p\nu }{\nu -p}$ if $\nu >p$ or $\beta =\sigma +1$ if $%
\nu =p$. There exists $C$ such that
\begin{equation}
\left\| u\right\| _{L^{p}\left( X,m\right) }\leq C\left\| u\right\|
_{W},\qquad \left\| u\right\| _{L^{\beta }\left( X,m\right) }\leq C\left\|
u\right\| _{W}.
\end{equation}
Choose $\varepsilon <\frac{1}{2C^{p}}$; then
\begin{equation*}
\int\limits_{X}F\left( u\right) m\left( dx\right) \leq \varepsilon \left[
\int\limits_{X}\left| u\right| ^{p}m\left( dx\right) +\int\limits_{X}\left|
u\right| ^{\beta }m\left( dx\right) \right] +C_{\varepsilon
}\int\limits_{X}\left| u\right| ^{\beta }m\left( dx\right)
\end{equation*}
\begin{equation}
=\varepsilon \left( \left\| u\right\| _{L^{p}\left( X,m\right) }^{p}+\left\|
u\right\| _{L^{\beta }\left( X,m\right) }^{\beta }\right) +C_{\varepsilon
}\left\| u\right\| _{L^{\beta }\left( X,m\right) }^{\beta }\leq \varepsilon
\left( C^{p}\left\| u\right\| _{W}^{p}+C^{\beta }\left\| u\right\|
_{W}^{\beta }\right) +C_{\varepsilon }C^{\beta }\left\| u\right\|
_{W}^{\beta }
\end{equation}
and
\begin{equation*}
\varphi \left( u\right) =\frac{1}{2}\left\| u\right\|
_{W}^{p}-\int\limits_{X}F\left( u\right) m\left( dx\right) \geq \left( \frac{%
1}{2}-\varepsilon C^{p}\right) \left\| u\right\| _{W}^{p}-C^{\beta }\left(
\varepsilon +C_{\varepsilon }\right) \left\| u\right\| _{W}^{\beta }
\end{equation*}
\begin{equation}
\geq \rho ^{p}-C^{\beta }\left( \varepsilon +C_{\varepsilon }\right) \rho
^{\beta }
\end{equation}
and the result follows from the last inequality. We consider now the case $%
\nu <2$. From the assumption we obtain that for every $\varepsilon >0$ there
exists a constant $\delta >0$ such that
\begin{equation}
F\left( t\right) \leq \varepsilon \left| t\right| ^{p}
\end{equation}
for $\left| t\right| \leq \delta $. We observe that there exists $C$ such
that
\begin{equation}
\left\| u\right\| _{L^{p}\left( X,m\right) }\leq C\left\| u\right\|
_{W},\qquad \left\| u\right\| _{L^{\infty }\left( X,m\right) }\leq C\left\|
u\right\| _{W}.
\end{equation}
Choosing $\left\| u\right\| _{W}=\rho =\frac{\delta }{C}$, we have $\left\|
u\right\| _{L^{\infty }\left( X,m\right) }\leq \delta $; then
\begin{equation}
\int\limits_{X}F\left( u\right) m\left( dx\right) \leq \varepsilon
\int\limits_{X}\left| u\right| ^{p}m\left( dx\right) =\varepsilon \left\|
u\right\| _{L^{p}\left( X,m\right) }^{p}\leq \varepsilon C^{p}\left\|
u\right\| _{W}^{p}
\end{equation}
and
\begin{equation}
\varphi \left( u\right) =\frac{1}{2}\left\| u\right\|
_{W}^{p}-\int\limits_{X}F\left( u\right) m\left( dx\right) \geq \left( \frac{%
1}{2}-\varepsilon C^{p}\right) \left\| u\right\| _{W}^{p}\geq \rho ^{p}.
\end{equation}
The result follows from the last inequality. Let us prove the existence of $%
u_{0}\in X\backslash B_{\rho }$ such that $\varphi \left( u\right) \leq 0$.
Let $u_{0}\in D\left[ a\right] $ be the potential of the ball $B\left(
0,1\right) $ with respect to the ball $B\left( 0,2\right) $. Then $u_{0}$ is
in $W$ and $\left\| u_{0}\right\| _{W}\geq am\left( B\left( 0,1\right)
\right) >\rho $; we recall that
\begin{equation}
F\left( u_{0}\left( x\right) \right) \geq m\left| u_{0}\left( x\right)
\right| ^{\mu }
\end{equation}
for $x\in B\left( 0,1\right) $. Let $\gamma >1$; we have $u_{0}\left(
x\right) =1$ on $B\left( 0,1\right) $, so
\begin{equation*}
\varphi \left( \gamma u_{0}\right) =\frac{1}{2}\gamma ^{p}\left\|
u_{0}\right\| _{W}^{p}-\int\limits_{X}F\left( \gamma u_{0}\right) m\left(
dx\right) \leq \frac{1}{2}\gamma ^{p}\left\| u_{0}\right\|
_{W}^{p}-\int\limits_{B\left( 0,1\right) }F\left( \gamma u_{0}\right)
m\left( dx\right)
\end{equation*}
\begin{equation}
\leq \frac{1}{2}\gamma ^{p}\left\| u_{0}\right\| _{W}^{p}-m\gamma ^{\mu
}\int\limits_{B\left( 0,1\right) }\left| u_{0}\right| ^{\mu }m\left(
dx\right) \leq \frac{1}{2}\gamma ^{p}\left\| u_{0}\right\| _{W}^{p}-m\gamma
^{\mu }m\left( B\left( 0,1\right) \right) .
\end{equation}
Since $\mu >p$ we have for $\gamma >\gamma _{0}$, $\gamma _{0}$ suitable, we
have $\varphi \left( \gamma u_{0}\right) <0$. The proof is completed with
the application of the Mountain Pass Theorem.
\end{proof}

\end{document}